\documentclass[12pt]{article}

\usepackage{amsthm,amsfonts,amsmath,amssymb}
\usepackage[dvips]{graphicx}
\usepackage{alltt,makeidx,enumerate}
\newtheorem{thm}{Theorem}[section]

\newtheorem{prop}[thm]{Proposition}

\theoremstyle{definition}

\newtheorem{remark}[thm]{Remark}

\newtheorem{example}[thm]{Example}
\numberwithin{equation}{section}
\newenvironment{prueba}{\noindent\textit{Proof}:}{  \qed\\\indent}

\newcommand{\me}{\mathrm{e}}

\newcommand{\dif}{\mathrm{d}}
\newcommand{\Dif}{\mathrm{D}}

\newcommand{\limite}[2]{\lim_{#1\rightarrow #2}}
\newcommand{\deru}[1]{#1^\prime}
\newcommand{\derd}[1]{#1^{\prime\prime}}

\newcommand{\abs}[1]{\left\vert#1\right\vert}
\newcommand{\set}[1]{\left\{#1\right\}}

\newcommand{\Real}{\mathbb R}

\newcommand{\To}{\rightarrow}

\newcommand{\ts}[1]{\mathbf{#1}}

\newcommand{\pd}[2]{\frac{\partial#1}{\partial#2}}

\makeindex

\begin{document}
\title{Chaotic dynamics of a three particle array under Lennard--Jones type
forces and a fixed area constraint}

\author{Pablo V. Negr\'on--Marrero\footnote{email: pablo.negron1@upr.edu}\\
Department of Mathematics\\University of Puerto Rico\\
Humacao, PR 00791-4300}
\date{}
\maketitle

\begin{abstract}
We consider the dynamical problem for a system of three particles in which the 
inter--particle forces are given as the gradient of a Lennard--Jones type 
potential. Furthermore we assume that the three particle array is subject to 
the constraint of fixed area. The corresponding mathematical problem is that of 
a conservative dynamical system over the manifold determined by the area 
constraint.
We study numerically the stability of this system. In particular, using the 
recently introduced measure of chaos by Hunt and Ott (2015), we study 
numerically the possibility of chaotic behaviour for this system.
\end{abstract}



\section{Introduction}
In this paper we consider the dynamical problem for a system of three 
particles, where the interactions between the particles of the system are due 
to forces given as the gradient of a Lennard--Jones \cite{Le1924} type 
potential. Furthermore we assume that the three particle array is subjected to 
the additional constraint that the area of the triangle generated by them is 
fixed. The motivation to consider the area constraint comes from the following 
phenomena observed both in laboratory experiments and molecular dynamics 
simulations (see e.g., \cite{BaNoSt2008,BlKa1975}). As the density of a fluid 
is progressively lowered (keeping the temperature constant), there is a certain 
``critical'' density such that if the density of the fluid is lower than this 
critical value, then bubbles or regions with very low density appear within the 
fluid. This phenomenon is usually called ``cavitation'' and it has been 
extensively studied as well in solids. (See for instance \cite{Fond2001, 
HoPo1995} for discussions and further references.)

The dynamical system under consideration in this paper is an example of a 
conservative dynamical system over a manifold (cf. \eqref{dyn_sist_full}). The 
numerical 
techniques to compute approximate solutions for these types of problems are 
based essentially on variations of a predictor--corrector (or projection) 
method. We refer to \cite{Ha2001}, \cite{Ha2011} for a discussion on the 
existence and uniqueness, dependence on initial data, as well as numerical 
schemes for general dynamical systems over manifolds. The assumptions that 
one particle is fixed at the origin and that another one moves along one of the 
coordinate axes, simplify greatly the area constraint in our problem. Using 
this we show in Section \ref{sec:4} that our original constrained dynamical 
system reduces to a standard (non-constrained) dynamical system 
(cf. \eqref{redsistp}) involving the area parameter. 

In Section \ref{sec:3} we give a characterization of the equilibrium points of 
\eqref{dyn_sist_full} or equivalently \eqref{redsistp}. We show that these 
equilibrium points coincide with those studied in \cite{NeMaLoSe2019}. In that 
paper the equilibrium problem is treated as a bifurcation problem with the 
specified area as the bifurcation parameter. They showed the existence of a 
family 
of equilibrium states corresponding to equilateral triangles, and that these 
equilibrium points are stable up to a certain critical area $A_c$ after which 
they become unstable. Furthermore, using techniques of equivariant bifurcation 
theory they showed that bifurcating from the point corresponding to $A_c$, there 
are solutions curves corresponding to isosceles triangles the stability of which 
could only be determined numerically. Moreover for an instance of the 
Lennard--Jones potential, they showed numerically that stable equilibrium points 
corresponding to scalene triangles exist for a range of values of the area 
parameter.

The study of \textit{chaotic} behaviour in dynamical systems is a very 
important and active area in this field. Verifying whether or not a particular 
system is chaotic can be a difficult task. In Section \ref{sec:chaos} we 
review a new criteria for chaos introduced recently by Hunt and Ott 
\cite{HuOt2015}. This new characterization of chaos is based on the so called 
\textit{expansion entropy} (cf. \eqref{expent}) which does not require the 
identification of a compact invariant set, and naturally leads itself to a 
computational method for detecting chaotic behaviour in a dynamical system. In 
Section \ref{sec:chaos} we describe this numerical scheme and give several 
examples of the use of this method applied to various dynamical systems of 
known 
chaotic or non--chaotic behaviour. Finally in Section \ref{sec:6} we use this 
method to show that the dynamical system for the three particle array under the 
area constraint, exhibits chaotic behaviour essentially for all values of the 
area parameter in the area constraint.

\section{Problem formulation}\label{sec:2}
We consider a system of three equal particles interacting via an
inter--particle potential $\phi$. The potentials we are concerned are those of 
the form:
\begin{equation}\label{LJ-potential}
\phi(r)=\dfrac{c_1}{r^{\delta_1}}-\dfrac{c_2}{r^{\delta_2}},
\end{equation}
where $c_1,c_2$ are positive constants and $\delta_1>\delta_2>2$. These 
constants determine the physical properties of the particles or molecules in the 
array. (The classical \textit{Lennard--Jones} \cite{Le1924} potential is 
obtained upon setting $\delta_1=12$ and $\delta_2=6$.) We study the planar 
dynamics of such a system subject to the constraint of fixed area for the 
triangular array. 

Let $\ts{r}_i(t)$, $i=1,2,3$ be the position vectors of the particles, and 
$\ts{r}_{ij}(t)=\ts{r}_i(t)-\ts{r}_j(t)$ for $i<j$. Let
\[
r_i=||\ts{r}_i||=(\ts{r}_i\cdot\ts{r}_i)^\frac{1}{2},\quad
r_{ij}=||\ts{r}_{ij}||.
\]
If $m$ is the particle mass, then the kinetic and potential energies of the
system are given respectively by:
\begin{equation}\label{kinpotener}
K=\frac{m}{2}\sum_{i=1}^3||\dot{\ts{r}}_i||^2,\quad 
U=\sum_{i<j}\phi(r_{ij}).
\end{equation}
Thus we seek to minimize the Lagrangian $L=K-U$ over the time interval $[0,t_f]$ which leads us to consider the functional:
\[
 E(\ts{r}_1,\ts{r}_2,\ts{r}_3)=\int_0^{t_f}(K-U)\,\dif t
\]

The square of the area of a triangle with sides $a,b,c$ is given according to Heron's formula by
\begin{equation}\label{Heron}
\Gamma(a,b,c)=\frac{1}{8}(a^2b^2+a^2c^2+b^2c^2)-\frac{1}{16}(a^4+b^4+c^4).
\end{equation}
For any $A>0$ we let 
\begin{equation}\label{areaconst}
g(\ts{r}_1,\ts{r}_2,\ts{r}_3)=\Gamma(r_{12},r_{13},r_{23})-A^2.
\end{equation}
Thus if we set $g=0$, then we would be specifying that the triangle defined by the positions of
the three particles is of fixed area $A$. We shall need the following formulas for the partial derivatives of $g$:
\begin{subequations}
\begin{eqnarray}
\pd{g}{\ts{r}_1}&=&\dfrac{1}{r_{12}}\,\pd{\Gamma}{a}(\ts{r}_1-\ts{r}_2)+
\dfrac{1}{r_{13}}\,\pd{\Gamma}{b}(\ts{r}_1-\ts{r}_3),
\label{pdg1}\\
\pd{g}{\ts{r}_2}&=&-\dfrac{1}{r_{12}}\,\pd{\Gamma}{a}(\ts{r}_1-\ts{r}_2)+
\dfrac{1}{r_{23}}\,\pd{\Gamma}{c}(\ts{r}_2-\ts{r}_3),
\label{pdg2}\\
\pd{g}{\ts{r}_3}&=&-\dfrac{1}{r_{13}}\,\pd{\Gamma}{b}(\ts{r}_1-\ts{r}_3)-
\dfrac{1}{r_{23}}\,\pd{\Gamma}{c}(\ts{r}_2-\ts{r}_3),
\label{pdg3}
\end{eqnarray}
\end{subequations}
where the partial derivatives of $\Gamma$ are evaluated at 
$(r_{12},r_{13},r_{23})$. Note that
\begin{equation}\label{symg}
 \pd{g}{\ts{r}_1}+\pd{g}{\ts{r}_2}+\pd{g}{\ts{r}_3}=\ts{0}.
\end{equation}
\subsection{The constrained problem}
We assume that the particle corresponding to $i=3$ is fixed at the origin,
while that for $i=2$ is restricted to move along the $x$--axis. Using the notation
$\ts{r}_i=(u_i,w_i)$, $i=1,2,3$, our variational problem becomes that of
finding $C^1$ functions $\ts{r}_1,\ts{r}_2,\ts{r}_3$ such that:
\[
\begin{array}{c}
\min~E(\ts{r}_1,\ts{r}_2,\ts{r}_3)\\
\mbox{subject~to~}g(\ts{r}_1,\ts{r}_2,\ts{r}_3)=0,\quad\ts{r}_3=\ts{0}
,\quad w_2=0.
\end{array}
\]
The Euler--Lagrange equations for the above variational problem are given by:
\begin{eqnarray}
m\ddot{\ts{r}}_1&=&-\dfrac{\deru{\phi}(r_{12})}{r_{12}}(\ts{r}_1-\ts{r}_2)-
\dfrac{\deru{\phi}(r_{13})}{r_{13}}(\ts{r}_1-\ts{r}_3)-\lambda\pd{g}{\ts{r}_1}
(\ts{r}_1,\ts{r}_2,\ts{r}_3),\nonumber\\
m\ddot{\ts{r}}_2&=&\dfrac{\deru{\phi}(r_{12})}{r_{12}}(\ts{r}_1-\ts{r}_2)-
\dfrac{\deru{\phi}(r_{23})}{r_{23}}(\ts{r}_2-\ts{r}_3)-\lambda\pd{g}{\ts{r}_2}
(\ts{r}_1,\ts{r}_2,\ts{r}_3)-\mu\ts{e}_2,\label{dyn_sist_full}\\
m\ddot{\ts{r}}_3&=&\dfrac{\deru{\phi}(r_{13})}{r_{13}}(\ts{r}_1-\ts{r}_3)+
\dfrac{\deru{\phi}(r_{23})}{r_{23}}(\ts{r}_2-\ts{r}_3)-\lambda\pd{g}{\ts{r}_3}
(\ts{r}_1,\ts{r}_2,\ts{r}_3)-\nu_1\ts{e}_1-\nu_2\ts{e}_2,\nonumber\\
0&=&g(\ts{r}_1,\ts{r}_2,\ts{r}_3),\quad\ts{r}_3=\ts{0},\quad w_2=0.\nonumber
\end{eqnarray}
Here $\lambda$, $\nu_1, \nu_2,\mu$ (which are functions of $t$) are the
Lagrange multipliers corresponding to the constraints in the problem, and
$\set{\ts{e}_1,\ts{e}_2}$ is the standard basis for $\Real^2$. From these
equations and \eqref{symg}, it follows now that
\[
m\sum_{i=1}^3\ddot{\ts{r}}_i(t)=-\nu_1\ts{e}_1-(\nu_2+\mu)\ts{e}_2,
\]
which describes the motion of the center of mass of the system. The system 
\eqref{dyn_sist_full} is conservative as the energy $K+U$ is conserved along 
solutions.

Applying the constraints $\ts{r}_3=\ts{0}$ and $w_2=0$, the equations above
simplify to:
\begin{eqnarray*}
m\ddot{\ts{r}}_1&=&-\dfrac{\deru{\phi}(r_{12})}{r_{12}}(\ts{r}_1-\ts{r}_2)-
\dfrac{\deru{\phi}(r_{1})}{r_{1}}\ts{r}_1-\lambda\pd{g}{\ts{r}_1},\\
m\ddot{u}_2&=&\dfrac{\deru{\phi}(r_{12})}{r_{12}}(u_1-u_2)-
\dfrac{\deru{\phi}(r_{2})}{r_{2}}u_2
-\lambda\pd{g}{u_2},\\
0&=&\dfrac{\deru{\phi}(r_{12})}{r_{12}}w_1
-\lambda\pd{g}{w_2}-\mu,\\
\ts{0}&=&\dfrac{\deru{\phi}(r_{1})}{r_{1}}\ts{r}_1+
\dfrac{\deru{\phi}(r_{2})}{r_{2}}\ts{r}_2-\lambda\pd{g}{\ts{r}_3}
-\nu_1\ts{e}_1-\nu_2\ts{e}_2,\\
0&=&g.
\end{eqnarray*}
where $g$ and all of its partial derivatives are evaluated at 
$(u_1,w_1,u_2,0,0,0)$. Note that since $r_2=|u_2|$, the first two equations 
together with the last one, can be solved independently of the third and fourth 
equations. Once $\ts{r}_1$, $u_2$, and $\lambda$ are determined, one can find
$\nu_1,\nu_2, \mu$ using the third and fourth equations. Thus we are lead to
consider the system in component form:
\begin{eqnarray}
m\ddot{u}_1&=&-\dfrac{\deru{\phi}(r_{12})}{r_{12}}(u_1-u_2)-
\dfrac{\deru{\phi}(r_{1})}{r_{1}}u_1-\lambda\pd{\gamma}{u_1}(u_1,w_1,u_2),
\nonumber\\
m\ddot{w}_1&=&-\left[\dfrac{\deru{\phi}(r_{12})}{r_{12}}+
\dfrac{\deru{\phi}(r_{1})}{r_{1}}\right]w_1-\lambda\pd{\gamma}{w_1}(u_1,w_1,
u_2),\label{redsist}\\
m\ddot{u}_2&=&\dfrac{\deru{\phi}(r_{12})}{r_{12}}(u_1-u_2)-
\dfrac{\deru{\phi}(r_{2})}{r_{2}}u_2-\lambda\pd{\gamma}{u_2}(u_1,w_1,u_2),
\nonumber\\
0&=&\gamma(u_1,w_1,u_2),\nonumber
\end{eqnarray}
where $\gamma(u_1,w_1,u_2)=g(u_1,w_1,u_2,0,0,0)$. 

We recall that in these equations
\begin{equation}\label{redrs}
 r_{12}=\sqrt{(u_1-u_2)^2+w_1^2},\quad r_1=\sqrt{u_1^2+w_1^2},\quad r_2=|u_2|.
\end{equation}
It is easy to check now that
\[
 \gamma(u_1,w_1,u_2)=\frac{1}{4}w_1^2u_2^2-A^2.
\]
Thus the area constraint $\gamma=0$ is simply the square of the familiar formula 
of ``area equals one half base times height''. We shall exploit this very 
shortly but before doing so, we pause to characterize the equilibrium points of 
the system \eqref{redsist}.

\section{The equilibrium points}\label{sec:3}
We now show that the equilibrium points of the system \eqref{redsist} 
are precisely those studied in \cite{NeMaLoSe2019}. We remark that while 
$\lambda$ in \eqref{redsist} is in general a function of $t$, in the following 
discussion it is a constant corresponding to a steady state.
\begin{prop}\label{eqredsist}
The equilibrium points of the system \eqref{redsist}
correspond to position vectors $\ts{r}_1=(u_1,w_1)$ and $\ts{r}_2=(u_2,0)$ for
which $a=r_{12}$, $b=r_1$, and $c=r_2$ (cf. \eqref{Heron}, \eqref{redrs})
satisfy:
\begin{equation}\label{equieq}
\left\{\begin{array}{rcl}
\deru{\phi}(a)+\lambda \Gamma_{,a}&=&0,\\
\deru{\phi}(b)+\lambda \Gamma_{,b}&=&0,\\
\deru{\phi}(c)+\lambda \Gamma_{,c}&=&0,\\
\Gamma(a,b,c)&=&A^2,
\end{array}\right.
\end{equation}
where $\Gamma_{,a}=\pd{\Gamma}{a}$, etc.
\end{prop}
\begin{prueba}
After setting the right hand side of \eqref{redsist} equal to zero, we are led 
to
\begin{subequations}\label{equil}
\begin{eqnarray}
\dfrac{\deru{\phi}(a)}{a}(u_1-u_2)+
\dfrac{\deru{\phi}(b)}{b}u_1+\lambda\pd{\gamma}{u_1}&=&0,\label{equi1}\\
\dfrac{\deru{\phi}(a)}{a}w_1+
\dfrac{\deru{\phi}(b)}{b}w_1+\lambda\pd{\gamma}{w_1}&=&0,\label{equi2}\\
-\dfrac{\deru{\phi}(a)}{a}(u_1-u_2)+
\dfrac{\deru{\phi}(c)}{c}u_2+\lambda\pd{\gamma}{u_2}&=&0,\label{equi3}\\
\gamma(u_1,w_1,u_2)&=&0,\label{equi4}
\end{eqnarray}
\end{subequations}
where we have set $a=r_{12}$, $b=r_1$, and $c=r_2$. The constraint 
$\gamma(u_1,w_1,u_2)=0$ is equivalent to
$\Gamma(a,b,c)=A^2$. Since $A>0$ is the area of the triangle determined by
$\ts{r}_1=(u_1,w_1)$, $\ts{r}_2=(u_2,0)$, and the origin, it follows that
\[
 w_1\ne0,\quad u_2\ne0.
\]
Now using that $w_2=0$ and $\ts{r}_3=\ts{0}$ together with the definitions of 
$a,b,c$, we get from \eqref{pdg1} and the first component of \eqref{pdg2} that:
\begin{eqnarray*}
\pd{\gamma}{u_1}&=&\dfrac{1}{a}\,\Gamma_{,a}(u_1-u_2)+
\dfrac{1}{b}\,\Gamma_{,b}u_1,\\
\pd{\gamma}{w_1}&=&\dfrac{1}{a}\,\Gamma_{,a}w_1+
\dfrac{1}{b}\,\Gamma_{,b}w_1,\\
\pd{\gamma}{u_2}&=&-\dfrac{1}{a}\,\Gamma_{,a}(u_1-u_2)+
\dfrac{1}{c}\,\Gamma_{,c}u_2.
\end{eqnarray*}
If we now substitute these equations into \eqref{equil}, and
add and subtract \eqref{equi1} and \eqref{equi3}, then the first three
equations of \eqref{equil} are equivalent to:
\begin{eqnarray}
H_bu_1+H_cu_2&=&0,\nonumber\\
(2H_a+H_b)u_1-(2H_a+H_c)u_2&=&0,\label{auxsist1}\\
(H_a+H_b)w_1&=&0,\nonumber
\end{eqnarray}
where
\[
 H_a=\frac{1}{a}\left(\deru{\phi}(a)+\lambda \Gamma_{,a}\right),\mbox{  etc.}
 \]
 Since $w_1\ne0$ we must have that $H_a=-H_b$. Since $u_2\ne0$ the determinant
of the coefficient matrix of the first two equations in the system
\eqref{auxsist1} must be zero. A calculation shows that this determinant is
$2H_a^2$. Hence $H_a=0$, and since $H_b=-H_a$, then $H_b=0$. It follows now
from the first equation of \eqref{auxsist1} and using that $u_2\ne0$, that
$H_c=0$. Since $H_a=H_b=H_c=0$ is equivalent to the first three equations in
\eqref{equieq}, the result follows.
\end{prueba}
For future reference we record here the basic result in \cite{NeMaLoSe2019}
concerning the existence and multiplicity of solutions of the system
\eqref{equieq}.
\begin{thm}\label{equi_states}
For any value of $A>0$, the system \eqref{equieq} has a solution of with
$a=b=c=a_A$ corresponding to an equilateral triangle with corresponding
Lagrange--multiplier $\lambda_A$ where:
\begin{equation}\label{symcpt}
a_A=\frac{2\sqrt{A}}{\sqrt[4]{3}},\quad\lambda_A=-\frac{4\deru{\phi}(a_A)}{
a_A^3 } .
\end{equation}
This equilibrium point is stable, that is, a minimizer of the potential energy
functional $U$ in \eqref{kinpotener} subject to the fixed area constraint, if
and only if
\begin{equation}\label{stab-cond}
\rho(A)\equiv\derd{\phi}(a_A)+\frac{3}{a_A}\deru{\phi}(a_A)>0.
\end{equation}
Moreover, if $A_0$ is a simple root of $\rho(\cdot)$, then there exist three
branches of solutions corresponding to isosceles triangles, bifurcating from the
branch of equilateral triangles at the point where $A=A_0$.
\end{thm}
The stability of the bifurcating branches in Theorem \ref{equi_states} can only
be determined numerically. In \cite{NeMaLoSe2019} they give numerical examples
in which these bifurcations are of the trans--critical type and also that there
can be secondary bifurcations into stable scalene triangles.

\section{The reduced problem}\label{sec:4}
The system \eqref{redsist} is an example of a dynamical system over a manifold. 
One can in principle solve this system directly using some of the numerical 
techniques for these types of problems, that essentially are based on variations 
of a predictor--corrector (or projection) method (cf. 
\cite{Ha2001}, \cite{Ha2011}). However 
because of the simplification of the constraint $\gamma=0$ in \eqref{redsist}, 
the system \eqref{redsist} can be reduced further to one in the variables 
$(u_1,w_1)$. This simplifies greatly the calculations in Section \ref{sec:6} 
when we apply to our problem the new criteria for detecting chaos in 
\cite{HuOt2015}.

As we mentioned before, the constraint $\gamma=0$ reduces to:
\[
 \frac{1}{4}w_1^2u_2^2-A^2=0,
\]
Since $w_1\ne0$ and $u_2\ne0$ and $A>0$, we may assume that $w_1>0$ and $u_2>0$. It follows now that
\begin{equation}\label{pdrg}
 \pd{\gamma}{u_1}=0,\quad\pd{\gamma}{w_1}=\frac{1}{2}w_1u_2^2,\quad 
\pd{\gamma}{u_2}=\frac{1}{2}w_1^2u_2.
\end{equation}
Using these and the fact that the constraint $\gamma=0$ is equivalent to 
$w_1u_2=2A$, we can eliminate $\lambda$ and the equation of $u_2$ from 
\eqref{redsist}, to get the following reduced system for $u_1,w_1$:
\begin{subequations}\label{redsistr1}
\begin{eqnarray}
m\ddot{u}_1&=&-\dfrac{\deru{\phi}(r_{12})}{r_{12}}(u_1-u_2)-
\dfrac{\deru{\phi}(r_{1})}{r_{1}}u_1,\label{redsistr1a}\\
m\left[w_1^2+u_2^2\right]\ddot{w}_1&=&2A \left[
-\dfrac{\deru{\phi}(r_{12})}{r_{12}}(u_1-u_2)+\deru{\phi}(u_2)\right]\nonumber\\
&&-\left[\dfrac{\deru{\phi}(r_{12})}{r_{12}}+\dfrac{\deru{\phi}(r_{1})}{r_{1}} \right]w_1^3+ 
8mA^2\dfrac{(\dot{w}_1)^2}{w_1^3}.\label{redsistr1b}
\end{eqnarray}
\end{subequations}
In these equations $w_1>0$ and any instance of $u_2$ should be replaced with $2A/w_1$. 

For the numerical calculations of the next section, we shall need to transform 
the system \eqref{redsistr1} into one of first order. If we let $v_1=\dot{u}_1$ 
and $v_2=\dot{w}_1$, then the above system is equivalent to:
\begin{subequations}\label{redsistp}
\begin{eqnarray}
\dot{u}_1&=&v_1,\quad\dot{w}_1=v_2,\nonumber\\
m\dot{v}_1&=&-\dfrac{\deru{\phi}(r_{12})}{r_{12}}(u_1-u_2)-
\dfrac{\deru{\phi}(r_{1})}{r_{1}}u_1,\label{redsistpa}\\
m\left[w_1^2+u_2^2\right]\dot{v}_2&=&2A \left[
-\dfrac{\deru{\phi}(r_{12})}{r_{12}}(u_1-u_2)+\deru{\phi}(u_2)\right]\nonumber\\
&&-\left[\dfrac{\deru{\phi}(r_{12})}{r_{12}}+\dfrac{\deru{\phi}(r_{1})}{r_{1}} \right]w_1^3+ 
8mA^2\dfrac{v_2^2}{w_1^3}.\label{redsistpb}
\end{eqnarray}
\end{subequations}
Given any $(p,q,\alpha,\beta)$ with $q>0$, it follows from the standard 
existence and uniqueness theorem for ode's that these equations have a unique 
solution satisfying
\[
 u_1(0)=p,\quad w_1(0)=q,\quad v_1(0)=\alpha,\quad v_2(0)=\beta.
\]
If we let $\ts{r}=(u_1,w_1,v_1,v_2)$ and $\ts{p}=(p,q,\alpha,\beta)$, then we 
shall write
$\ts{r}(t;\ts{p})$ to denote the dependence of the solution of \eqref{redsistp}
on the initial conditions $\ts{p}$. It follows from standard results on the
dependence on initial values for the solutions of initial value problems (cf.
\cite{Ha1982}), that $\ts{r}(t;\ts{p})$ is a differentiable function of
$\ts{p}$.

\section{A criteria for chaos}\label{sec:chaos}
The study of \textit{chaotic} behaviour in dynamical systems is a very 
important and active
area in this field. A commonly used definition of chaos, originally formulated 
by \textit{Robert L. Devaney}, says that for a dynamical system to be 
classified as \textit{chaotic}, it must have the following properties 
\cite{HaKa2003}:
\begin{itemize}
\item
it must be sensitive to initial conditions
\item
it must be topologically mixing
\item
it must have dense periodic orbits
\end{itemize}
Verifying whether these properties hold or not for a particular system can be 
a difficult task. Recently Hunt and Ott \cite{HuOt2015} developed a new 
criteria for chaos based on the so called \textit{expansion entropy} (cf. 
\eqref{expent}) which does not require the identification of a compact 
invariant set. In this section we introduce some of the notions and definitions 
in \cite{HuOt2015} leading to a practical, from the numerical point of view, 
characterization of chaos for a dynamical system. 

We consider an autonomous dynamical system of the form:
\begin{equation}\label{adynsist}
 \dot{\ts{r}}(t)=\ts{f}(\ts{r}(t)),
\end{equation}
where $\ts{f}:\Omega\To\Real^n$ is a smooth function, and 
$\Omega\subset\Real^n$ is open with a smooth boundary. We denote by 
$\ts{r}(\cdot;\ts{p})$ the solution of \eqref{adynsist} such that 
$\ts{r}(0)=\ts{p}$. $\ts{r}(\cdot;\ts{p})$ is called the \textit{flow} or 
\textit{evolution operator} of \eqref{adynsist}. Let
\begin{equation}\label{derIC}
 \ts{u}(\cdot;\ts{p})=\Dif_{\ts{p}}\ts{r}(\cdot;\ts{p}).
\end{equation}
Note that $\ts{u}(t;\ts{p})$ is an $n\times n$ matrix for each $t$. Now since 
$\ts{r}(0;\ts{p})=\ts{p}$, it follows that $\ts{u}(0;\ts{p})=\ts{I}$, where 
$\ts{I}$ is the $n\times n$ identity matrix. After differentiating in 
\eqref{adynsist} with respect to $\ts{p}$, we get that
\[
 \dot{\ts{u}}(t)=\Dif_{\ts{r}}\ts{f}(\ts{r}(t))\ts{u}(t).
\]
Thus, the flow $\ts{r}(\cdot;\ts{p})$ and the matrix function 
$\ts{u}(\cdot;\ts{p})$ are solutions of the initial value problem:
\begin{subequations}\label{ruIC}
 \begin{eqnarray}
  \dot{\ts{u}}(t)&=&\Dif_{\ts{r}}\ts{f}(\ts{r}(t))\ts{u}(t),\label{ruIC1}\\
  \dot{\ts{r}}(t)&=&\ts{f}(\ts{r}(t)),\label{ruIC2}\\
  \ts{u}(0)&=&\ts{I},\quad\ts{r}(0)=\ts{p}.\label{ruIC3}
 \end{eqnarray}
\end{subequations}

The function $\ts{u}$ is central in the definition of chaos in \cite{HuOt2015}. 
In particular, let $G(\ts{u})$ be the product of the singular values of $\ts{u}$ 
that are greater than $1$. (If none of the singular values are greater than $1$, 
we set $G(\ts{u})=1$.) For any subset $S$ of $\Omega$, called a 
\textit{restraining region}, let
\[
 S_T=\set{\ts{p}\in S\,:\,\ts{r}(t;\ts{p})\in S,\quad t\in[0,T]},
\]
and define
\begin{equation}\label{orbgwrt}
 E_T(\ts{r},S)=\frac{1}{|S|}\int_{S_T}G(\ts{u}(T;\ts{p}))\,\dif\ts{p},
\end{equation}
where $|S|$ is the volume of $S$. The \textit{expansion entropy} 
$H_0(\ts{r},S)$ of the flow $\ts{r}$ over the set $S$ is defined as:
\begin{equation}\label{expent}
 H_0(\ts{r},S)=\limite{T}{\infty}\dfrac{\ln E_T(\ts{r},S)}{T}.
\end{equation}
According to \cite{HuOt2015} the system \eqref{adynsist} is \textit{chaotic} if 
for some restraining region $S$, we have that $H_0(\ts{r},S)>0$.
\begin{remark}\label{monoT}
The expansion entropy function has several nice properties, one of them being 
that $H_0(\ts{r},S')\le H_0(\ts{r},S)$ whenever $S'\subset S$. Thus if chaos is 
detected for any region $S'$, it will also be detected for any other region 
containing $S'$.
\end{remark}
\begin{remark}
The definitions in this section can be extended to non--autonomous systems and 
for any manifold $\Omega$ in $\Real^n$ (cf. \cite{HuOt2015}).
\end{remark}

\begin{example}
 We consider the constant coefficient linear system
 \begin{equation}\label{lin_sist}
   \dot{\ts{r}}(t)=A\ts{r}(t),
 \end{equation}
where $A$ is an $n\times n$ matrix. The solution of this system is given by
\[
 \ts{r}(t)=\me^{tA}\ts{p},\quad t\in\Real,
\]
where $\ts{p}\in\Real^n$ is arbitrary. Here
\[
 \me^{tA}=\sum_{k=0}^\infty\frac{t^k}{k!}\,A^k.
\]
It follows from the above representation of the solution $\ts{r}(\cdot)$ that 
the linear dynamical system above is non chaotic. We show now that the 
definition of chaos in \cite{HuOt2015} gives indeed a non chaotic behaviour in 
this case.

For simplicity of exposition we assume that $A$ is symmetric with eigenvalues 
$\lambda_1\le\lambda_2\le\cdots\le\lambda_n$ and corresponding orthonormal 
eigenvectors $\set{\ts{x}_1,\ts{x}_2,\ldots,\ts{x}_n}$. In this case, the 
matrix $\me^{tA}$ is symmetric as well with eigenvalues $\me^{\lambda_1t}$, 
$\me^{\lambda_2t}$,...,$\me^{\lambda_nt}$ with the same corresponding 
eigenvectors as $A$. Since the eigenvalues of $\me^{tA}$ are all positive, it 
follows that the singular values of $\me^{tA}$ coincide with its eigenvalues. 
Since 
the matrix $\ts{u}$ in \eqref{derIC} is given by $\me^{tA}$, we get that
\[
 G(\ts{u}(t;\ts{p}))=\me^{t\sum_{\lambda_i>0}\lambda_i},
\]
or one if the sum in the exponent is empty. In particular $G(\ts{u}(t;\ts{p}))$ 
is independent of $\ts{p}$.

Writing $\ts{p}=\sum_{i=1}^nc_i\ts{x}_i$, we get that
\[
 \ts{r}(t)=\sum_{i=1}^nc_i\me^{\lambda_it}\ts{x}_i. 
\]
For the restraining region, we consider the set
\[
S=\set{\ts{p}=\sum_{i=1}^nc_i\ts{x}_i\,:\,\abs{c_i}\le a},
\]
for any $a>0$. If all $\lambda_i$ are nonpositive, then $S_T=S$, 
$G(\ts{u}(T;\ts{p}))=1$, and we get that $E_T(\ts{r},S)=1$, and thus that 
$H_0(\ts{r},S)=0$.

Assume now that there are $p>0$ positive eigenvalues $\lambda_i$. In this case
\[
 S_T=\set{\sum_{i=1}^nc_i\ts{x}_i\,:\,\abs{c_i}\le a,\,1\le i\le n-p,\quad 
\abs{c_i}\le a\me^{-\lambda_iT},\,n-p+1\le i\le n},
\]
and $G(\ts{u}(T;\ts{p}))=\me^{T\sum_{i=n-p+1}^n\lambda_i}$. As 
$\abs{S_T}=(2a)^{n} \me^{-T\sum_{i=n-p+1}^n\lambda_i}$, we get that
\begin{eqnarray*}
 E_T(\ts{r},S)&=&\frac{1}{|S|}\int_{S_T}G(\ts{u}(T;\ts{p}))\,\dif\ts{p},\\
 &=&\frac{1}{(2a)^n}G(\ts{u}(T;\ts{p}))\abs{S_T}=1.
\end{eqnarray*}
Once again we get that that $H_0(\ts{r},S)=0$.

For an arbitrary compact restraining region $S^\prime$, we take $a$ in the 
definition of $S$ above sufficiently large such that $S^\prime\subset S$. It 
follows from Remark \ref{monoT} that $H_0(\ts{r},S')\le H_0(\ts{r},S)=0$, and 
thus that the linear dynamical system is non chaotic according to 
\cite{HuOt2015}.
\end{example}

The definition of $H_0(\ts{r},S)$ given in this section leads naturally to a 
very practical numerical scheme for estimating this quantity and thus to 
numerically test for chaos for a given dynamical system. If we take a set 
$\set{\ts{p}_1,\ldots,\ts{p}_N}$ of uniformly random vectors in $S$, then we can 
approximate \eqref{orbgwrt} by $\hat{E}_T(\ts{r},S)$ were
\begin{equation}\label{approxH0}
 \hat{E}_T(\ts{r},S)=\dfrac{1}{N}\sum_{{k=1} \atop {\ts{p}_k\in S_T}}^N 
G(\ts{u}(T;\ts{p}_k)).
\end{equation}
This sum is computed for different values of $T\in[0,T_{\mbox{max}}]$ for some 
prescribed $T_{\mbox{max}}$. The whole computation is repeated for different 
random samples of the $\set{\ts{p}_1,\ldots,\ts{p}_N}$. We then compute the 
average over the samples of $\ln[\hat{E}_T(\ts{r},S)]$, for each of the chosen 
$T$'s in $[0,T_{\mbox{max}}]$. From a plot of these averages vs $T$, we can 
estimate $H_0(\ts{r},S)$ as the asymptotic slope of this graph. (See 
\cite{HuOt2015}.) We close this section with an application of this numerical 
scheme to several dynamical systems whose possible chaotic or non--chaotic 
behaviour  
is known. The 
computations in this section and the rest of the paper were performed using the 
\texttt{ode45} routine of MATLAB using an event subroutine to detect when an 
orbit may exit $S$ for the first time.
\begin{example}
 As a first example we consider the special case of \eqref{lin_sist} in which
 \[
  A=\left[\begin{array}{rr}-0.45 &-0.55\\-0.55 &-0.45\end{array}\right].
 \]
We already shown that this system is non--chaotic according to the Hunt and 
Ott criteria. The coefficient matrix has eigenvalues 
$-1$ and $0.1$. The result of using the numerical scheme described above with 
$N=5000$ and $40$ samples, and 
restraining region
\[
 S=\set{(x,y)\,:\,\abs{x}\le10,\quad\abs{y}\le10},
\]
is shown in Figure \ref{fig:11}. The slope of approximately $8\times10^{-4}$ of 
the best line in this case is consistent with a non chaotic system.
\begin{figure}
\begin{center}
\scalebox{0.5}{\includegraphics{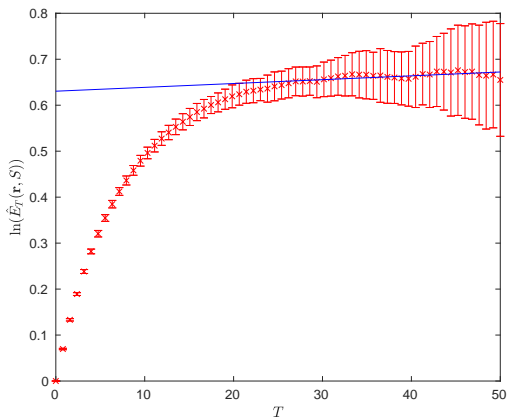}}
\end{center}
\caption{A simulation of the Hunt and Ott algorithm on a linear 
system with constant coefficients.}\label{fig:11}
\end{figure}
\end{example}

\begin{example}
 We consider the following dynamical system called a ``Sprott system of case 
A'':
\begin{eqnarray*}
 \dot{x}(t)&=&y,\\
 \dot{y}(t)&=&-x+yz,\\
 \dot{z}(t)&=&1-y^2.
\end{eqnarray*}
This system is an example of a chaotic dynamical system with no equilibrium 
points (cf. \cite{JaSpNa2015}, \cite{Sp1994}). We tested on this system the 
numerical scheme described above with $N=10000$ and $40$ samples, and 
restraining region
\[
 S=\set{(x,y,z)\,:\,\abs{x}\le10,\quad\abs{y}\le10,\quad\abs{z}\le10}.
\]
We show in Figure \ref{fig:8} a plot of the average  
$\ln[\hat{E}_T(\ts{r},S)]$ vs $T$ with their corresponding error bars (in terms 
of their variances). The positive slope of the best linear fitting 
in the picture (indicative of the chaotic behaviour of the system), is 
approximately $0.0252$. This number is an approximation of the sum of positive 
Lyapunov exponents over trajectories that remain in $S$ for all forward time 
\cite{HuOt2015}. 
\begin{figure}
\begin{center}
\scalebox{0.5}{\includegraphics{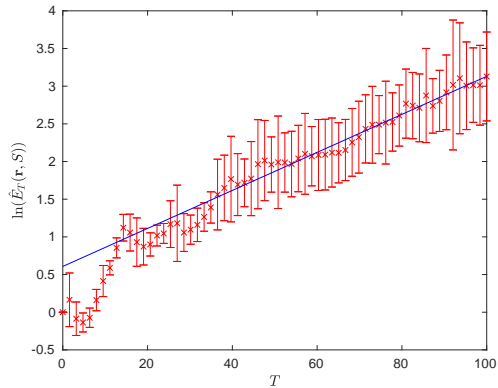}}
\end{center}
\caption{A simulation of the Hunt and Ott algorithm on a chaotic system 
with no equilibrium points.}\label{fig:8}
\end{figure}
\end{example}
\begin{example}
We consider now the motion of a double pendulum consisting of two point masses 
$m_1$, $m_2$ joined together by two massless shafts or bars of lengths $l_1$, 
$l_2$. (See Figure \ref{fig:9}.) We assume that there is no friction on the 
joints. If we let $\theta_i(t)$ to be the angle that the shaft $l_i$ makes with 
the vertical direction, $i=1,2$, then the equations for the motion of such a 
system are given after simplification by:
\begin{eqnarray*}
l_1\ddot{\theta}_1(t)&=&\frac{T_2}{m_1}\,
\sin(\theta_2(t)-\theta_1(t))-g\sin\theta_1(t),\\
l_2\ddot{\theta}_2(t)&=&-\frac{T_1}{m_1}\,
\sin(\theta_2(t)-\theta_1(t)),
\end{eqnarray*}
where $g$ is the constant of gravity acceleration and $T_1,T_2$, the tensions 
in the shafts, are solutions of the system
\[
\left[\begin{array}{cc}\frac{1}{m_1}&-\frac{1}{m_1}\,
\cos(\theta_2(t)-\theta_1(t))\\
-\frac{1}{m_1}\,\cos(\theta_2(t)-\theta_1(t))&\frac{1}{m_1}+\frac{1}{m_2}
\end{array}\right]\left[\begin{array}{c}T_1\\T_2\end{array}\right]=
\left[\begin{array}{c}l_1\dot{\theta}_1^2(t)+g\cos\theta_1(t)\\
       l_2\dot{\theta}_2^2(t)
      \end{array}\right].
\]
\begin{figure}
\begin{center}
\scalebox{0.5}{\includegraphics{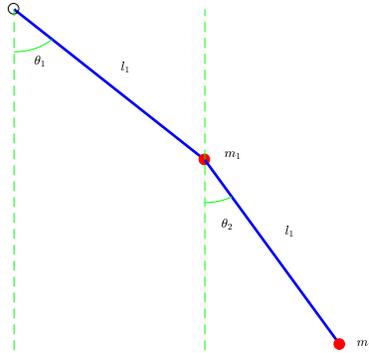}}
\end{center}
\caption{Frictionless double pendulum.}\label{fig:9}
\end{figure}
This system is well known to be chaotic (cf. \cite{De2011}, 
\cite{ShGrWiYo1992}). The Hunt and Ott algorithm was tested in this case with 
$N=1000$ and $20$ 
samples, and 
restraining region
\[ 
S=\set{(\theta_1,\dot{\theta}_1,\theta_2,\dot{\theta}_2)\,:\,\abs{\theta_i}
\le100 , \quad\abs {\dot{\theta}_i} \le20 , \quad i=1,2 } .
\]
In Figure \ref{fig:10} we show the results for this simulation. The slope of 
approximately $5.7016$ in this picture is consistent with the chaotic character 
of 
this system.
\begin{figure}
\begin{center}
\scalebox{0.5}{\includegraphics{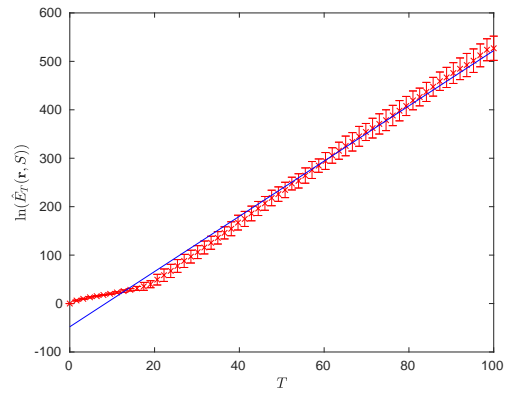}}
\end{center}
\caption{Simulation using the Hunt and Ott Frictionless double 
pendulum.}\label{fig:10}
\end{figure}
\end{example}

\section{Numerical results}\label{sec:6}
In this section we present various numerical results for the behaviour of 
solutions of \eqref{redsistp} when the potential $\phi$ is given by 
\eqref{LJ-potential}. In particular, using the Hunt and Ott method, we test the 
system \eqref{redsistp} for chaotic behaviour.

\begin{figure}
\begin{center}
\scalebox{0.5}{\includegraphics{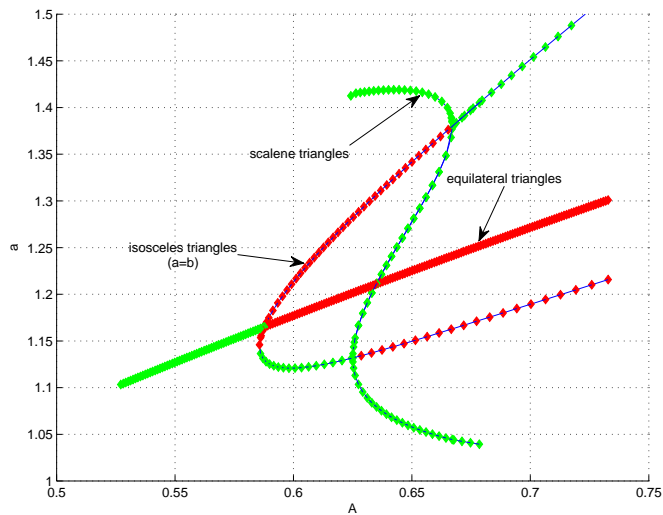}}
\end{center}
\caption{Bifurcation diagram for the equilibrium points of the system 
\eqref{redsistp} for the
Lennard-Jones potential \eqref{LJ-potential} with values 
\eqref{EXLJ}, green representing stable equilibria while red 
unstable ones.}\label{fig:1a}
\end{figure}

We consider the case of \eqref{LJ-potential} in which the parameters are given 
by 
\begin{equation}\label{EXLJ}
c_1=1,\quad c_2=2,\quad \delta_1=12,\quad  \delta_2=6. 
\end{equation}
The stability and multiplicity of the equilibrium states of \eqref{redsistp} 
(cf. 
Proposition \ref{eqredsist}) was fully analysed in \cite{NeMaLoSe2019} for 
this case. We briefly review those results that are more relevant to the present 
discussion. In Figure \ref{fig:1a} (reproduced from \cite{NeMaLoSe2019}) we show 
a projection of the bifurcation diagram of equilibrium states for the 
particular case \eqref{EXLJ}. (The notation here is as in Proposition 
\ref{eqredsist} where $a,b,c$ represent the sides of the triangular 
configuration of the array.) In this figure there are three bifurcation points 
at the following approximate values of the area parameter $A$ in 
\eqref{redsistp}: $A_0=0.5877$ (primary bifurcation) and $A_1=0.6251$, 
$A_2=0.6670$ (secondary bifurcations). The is also a turning point at 
$A_{-1}=0.5855$ corresponding to a trans--critical type bifurcation from $A_0$. 
The stability, multiplicity and type of the equilibrium point is summarized in 
Table \ref{tab:1}. The ``type'' column refers as to whether the shape of 
triangular array is an equilateral, isosceles, or scalene triangle. 

\begin{table}
\begin{center}
\begin{tabular}{|c|c|c|c|}\hline
$A$ & type & stability & multiplicity\\\hline\hline
$A\in(0,A_{-1})$ & equilateral & stable & one\\\hline\hline
 & equilateral & stable & one\\\cline{2-4}
\raisebox{0.0ex}[0pt]{$A\in(A_{-1},A_0)$} & isosceles & stable & 
three\\\cline{2-4}
 & isosceles & unstable & three\\\hline\hline
 & equilateral & unstable & one\\\cline{2-4}
\raisebox{0.0ex}[0pt]{$A\in(A_0,A_1)$} & isosceles & stable & 
three\\\cline{2-4}
 & isosceles & unstable & three\\\hline\hline
 & equilateral & unstable & one\\\cline{2-4}
\raisebox{0.0ex}[0pt]{$A\in(A_1,A_2)$} & scalene & stable & six\\\cline{2-4}
& isosceles & unstable & six\\\hline\hline
 & equilateral & unstable & one\\\cline{2-4}
\raisebox{0.0ex}[0pt]{$A\in(A_2,\infty)$} & isosceles & stable & 
three\\\cline{2-4}
 & isosceles & unstable & three\\\hline
\end{tabular}
\end{center}
\caption{Type, stability, and multiplicity for  the equilibrium points of 
\eqref{redsistp} for the values in \eqref{EXLJ}.}\label{tab:1}
\end{table}

\begin{figure}
\begin{center}
\scalebox{0.3}{\includegraphics{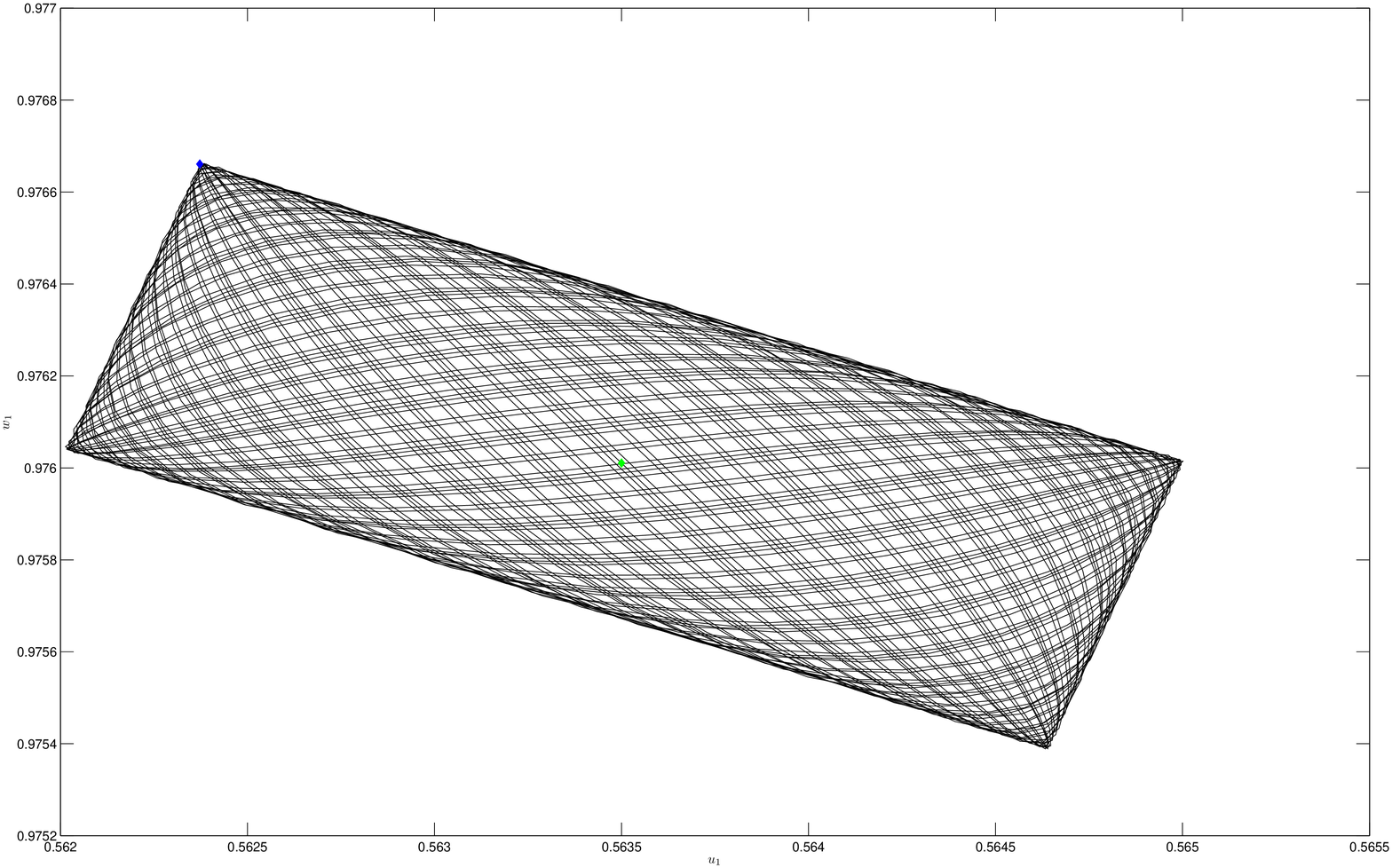}}
\end{center}
\caption{Sample orbit for initial point (blue) near the stable equilibrium point
(green) corresponding to $A=0.55$.}\label{fig:1}
\end{figure}

\begin{figure}
\begin{center}
\scalebox{0.3}{\includegraphics{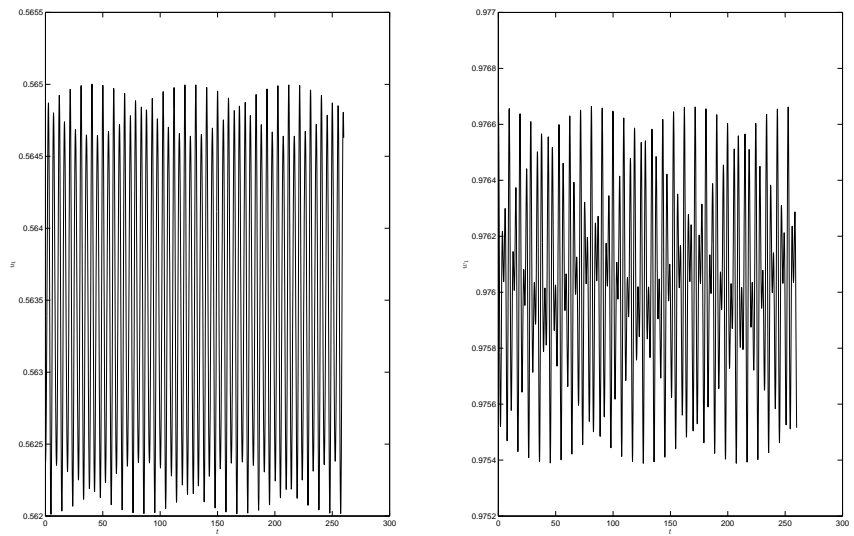}}
\end{center}
\caption{Graphs of $u_1$ and $w_1$ vs $t$ for the components of the orbit
in Figure \ref{fig:1} corresponding to $A=0.55$.}\label{fig:2}
\end{figure}

In the first simulation we examine the  dynamics of the system \eqref{redsistp} 
near the equilibrium point for the value $A=0.55$. This equilibrium point 
corresponds to an equilateral configuration of the array with sides $1.1270$ 
approximately. The corresponding equilibrium point for \eqref{redsistp} is 
$u_1=0.5635$, $w_1=0.9760$. We introduced a perturbation of $10^{-3}$ on one of 
the sides of the triangular array and used the resulting values of $u_1,w_1$, 
with velocities $v_1,v_2$ both $O(10^{-4})$, as initial conditions for 
\eqref{redsistp}. In Figure \ref{fig:1} we show the projection  onto the 
$u_1$--$w_1$ plane of the computed orbit of the dynamical system. Figure 
\ref{fig:2} shows the evolution of $u_1(t)$ and $w_1(t)$ as functions of time. 
The initial point is marked in blue and the equilibrium point in green. Both 
figures are consistent with a stable (not asymptotically stable) fixed point. 
Similar results are obtained for initial values close to other stable 
equilibrium corresponding to different values of $A$.

\begin{figure}
\begin{center}
\scalebox{0.3}{\includegraphics{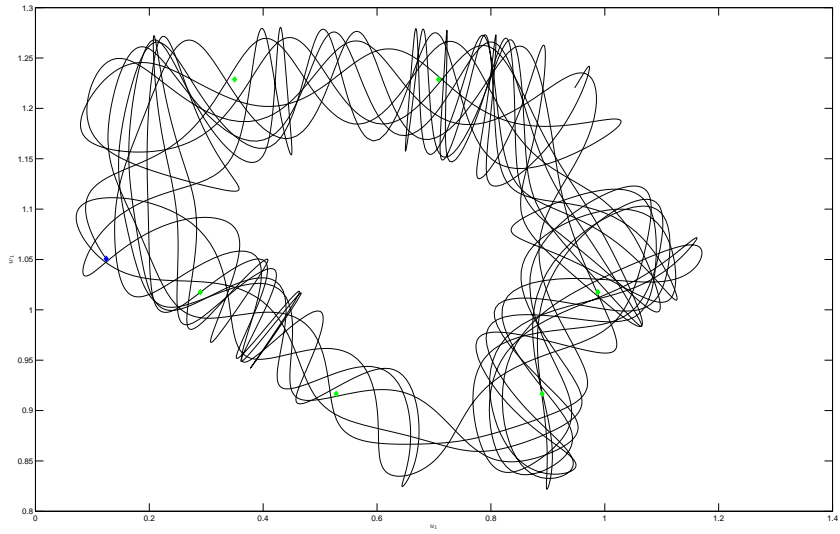}}
\end{center}
\caption{Sample orbit visiting all six stable equilibrium point
(green) corresponding to $A=0.65$.}\label{fig:3}
\end{figure}

\begin{figure}
\begin{center}
\scalebox{0.3}{\includegraphics{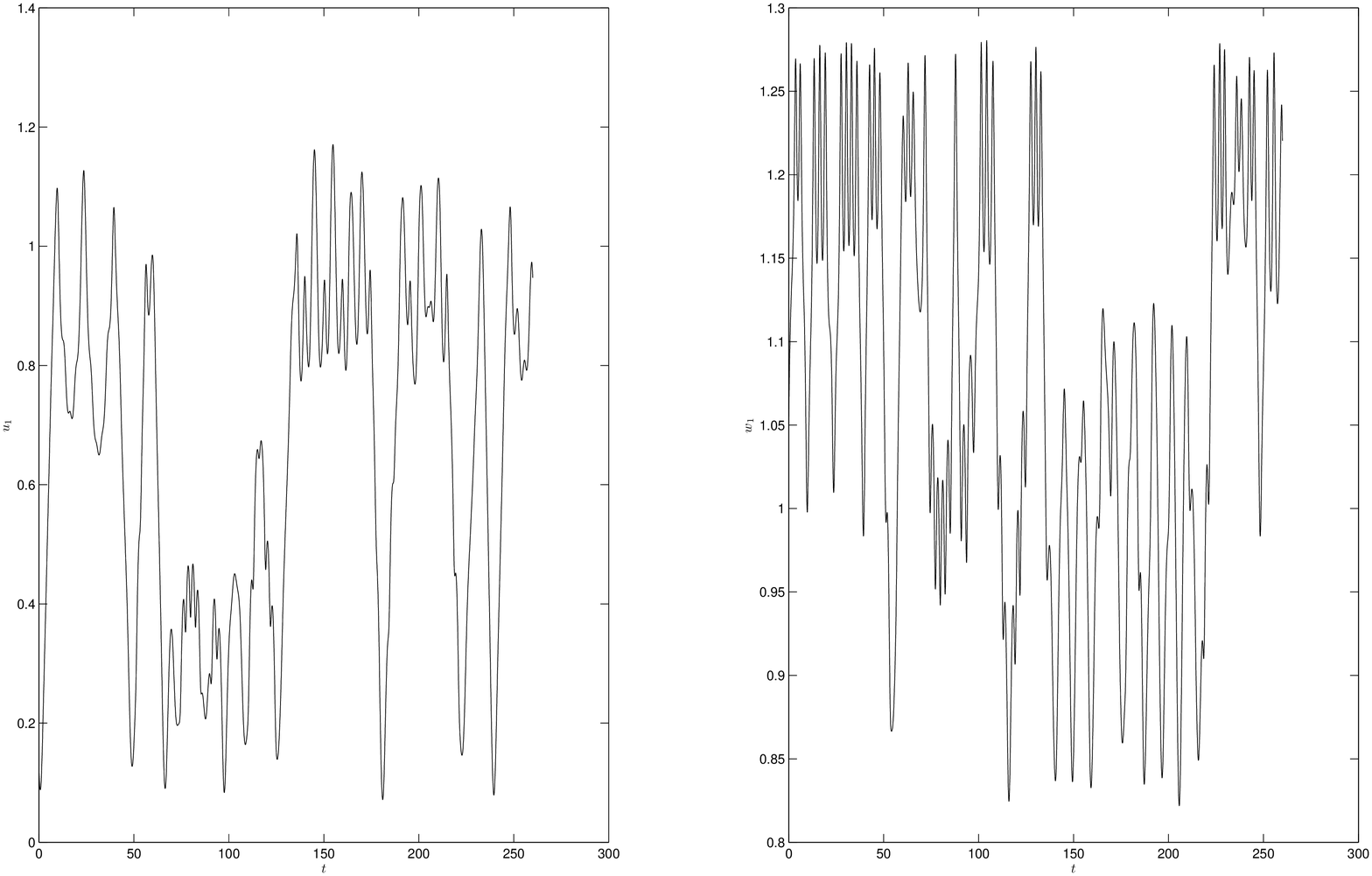}}
\end{center}
\caption{Graphs of $u_1$ and $w_1$ vs $t$ for the components of the orbit
in Figure \ref{fig:3} corresponding to $A=0.65$.}\label{fig:4}
\end{figure}

\begin{figure}
\begin{center}
\scalebox{0.3}{\includegraphics{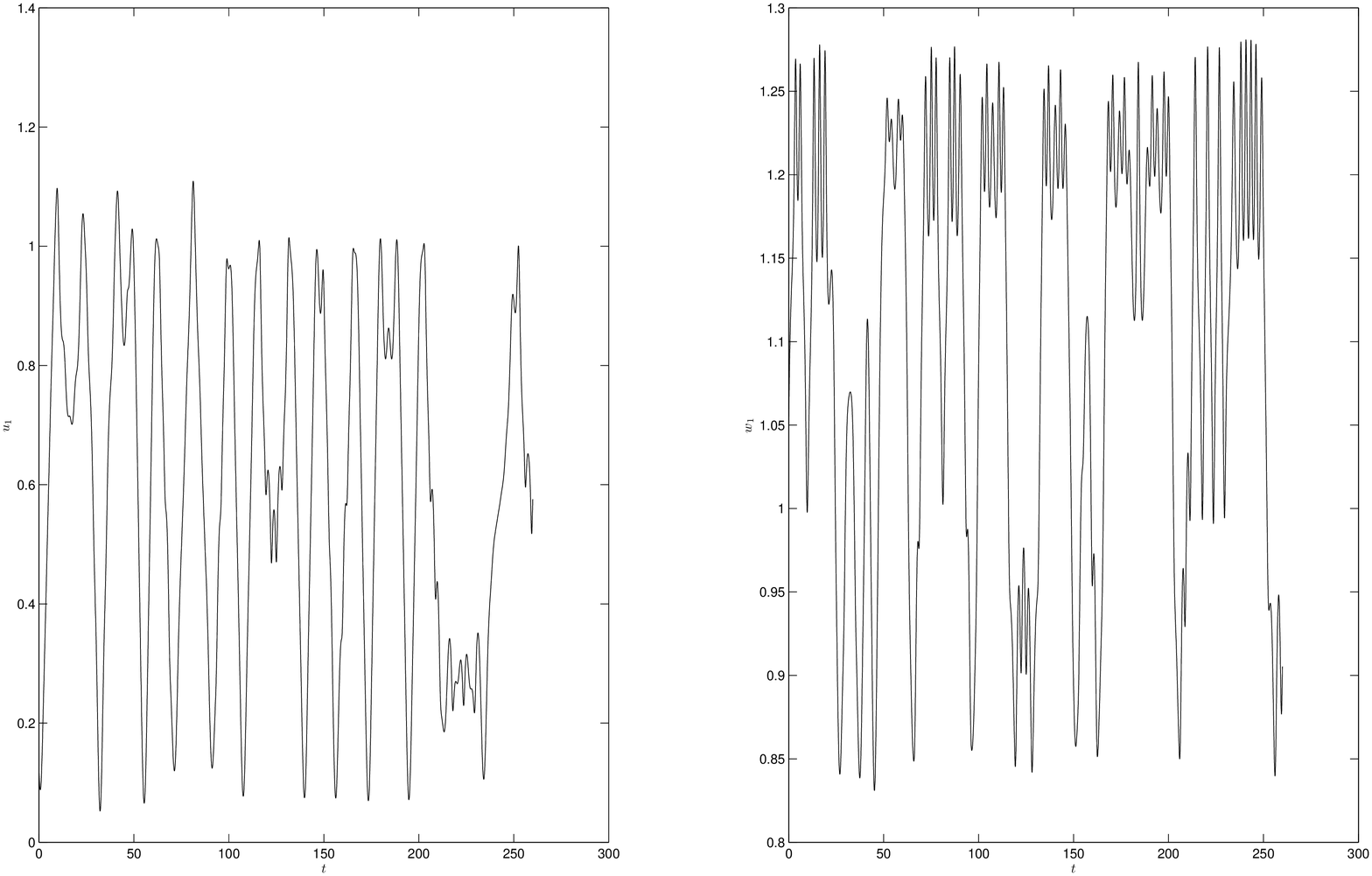}}
\end{center}
\caption{Graphs of $u_1$ and $w_1$ vs $t$ for the components of the orbit with
initial conditions as in Figure \ref{fig:4} plus a perturbation $O(10^{-4})$
corresponding to $A=0.65$.}\label{fig:5}
\end{figure}

We now consider the case in which $A=0.65$.  According to Table \ref{tab:1} we 
have one equilateral configuration which is unstable, six isosceles unstable 
configurations, and six stable equilibrium points 
corresponding to scalene triangles. For the scalene case, the triangle sides are 
$1.2776, 
1.4182, 1.0580$, the six triangular configurations obtained by permuting these 
numbers. In Figure \ref{fig:3} we show the projection onto the $u_1$--$w_1$ 
plane of an orbit generated with an initial point (indicated in blue) not 
necessarily close to any of the six stable equilibrium points (in green). Not 
that the 
orbit appears to visit ``regularly'' all the stable equilibrium points. In 
Figure 
\ref{fig:4} we show the evolution of $u_1(t)$ and $w_1(t)$ as functions of time, 
with the apparent random nature characteristic of a chaotic system. In Figure 
\ref{fig:5} we tested the orbit in Figure \ref{fig:4} for sensitivity to initial 
conditions. This figure was generated introducing a perturbation\footnote{Both 
orbits in Figures \ref{fig:4} and \ref{fig:5} were computed setting the 
absolute and relative error tolerances in ode45 to $10^{-10}$.} of $O(10^{-4})$ 
into  the initial condition used for Figure \ref{fig:4}. The resulting figure 
differs substantially from that in Figure \ref{fig:5}, again characteristic of a 
chaotic system.

\begin{figure}
\begin{center}
\scalebox{0.5}{\includegraphics{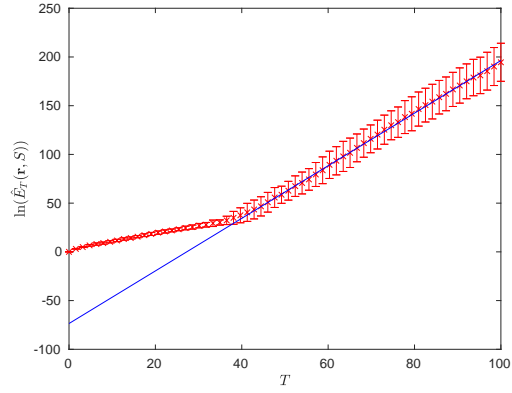}}
\end{center}
\caption{Graph of $\ln(\hat{E}_T)$ vs $T$ for the system \eqref{redsistp}
corresponding to $A=0.65$.}\label{fig:6}
\end{figure}

\begin{figure}
\begin{center}
\scalebox{0.5}{\includegraphics{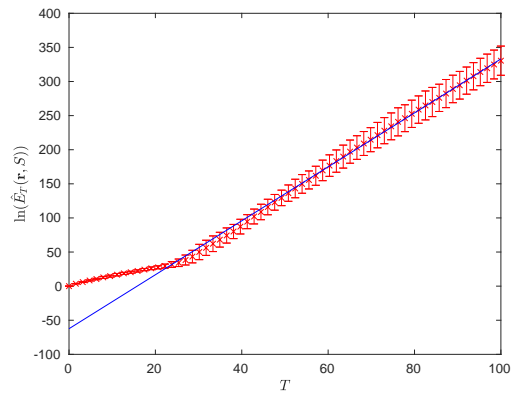}}
\end{center}
\caption{Graph of $\ln(\hat{E}_T)$ vs $T$ for the system \eqref{redsistp}
corresponding to $A=0.55$.}\label{fig:7}
\end{figure}

To test for possible chaotic behaviour of the system \eqref{redsistp}, we used 
the numerical scheme described at the end of Section \ref{sec:chaos} to 
approximate the expansion entropy $H_0$ for our system. Note that the entropy 
depends now on the area parameter $A$. As a restraining region we used:
\[
 S=\set{(u_1,w_1,v_1,v_2)\,:\,-2\le u_1,v_1,v_2\le 2,\quad 0<w_1\le2}.
\]
We computed \eqref{approxH0} with $N=5000$  and $50$ data samples. In Figure 
\ref{fig:6} we show a plot of the averages of the computed 
$\ln[\hat{E}_T(\ts{r},S)]$ for a range of values of $T$ corresponding to 
$A=0.65$. The red bars in the graph give intervals of plus or minus one sample 
standard deviation from each computed mean. As can be seen from the figure, the 
computed averages lie approximately on a line for large values of $T$. The slope 
of this line gives an approximation of $H_0(\ts{r},S)$ for our system. This 
slope is approximately $2.6988$ and thus according to the criteria in 
\cite{HuOt2015}, the system \eqref{redsistp} is chaotic for $A=0.65$. We 
performed a similar calculation for $A=0.55$. We show in Figure \ref{fig:7} the 
corresponding graph, again with the computed averages lying approximately on a 
line for large values of $T$, with slope of approximately $3.9548$ in this 
case. Thus the system is chaotic for this value of $A$ as well.

\section{Final comments}
The criteria for chaos given by Hunt and Ott \cite{HuOt2015} leads itself to a 
practical numerical method for detecting chaos in dynamical systems. The method 
is applicable to continuous as well as to discrete dynamical systems, even 
non--autonomous systems. The computations can be very intensive as the method 
requires a large number of random points over the restraining region, and this 
must be repeated another number of times in order to compute the required 
averages. However the method is easy to run in parallel which helps to 
reduce the computational time.

The values of $A=0.55$ and $0.65$ are typical values for different ranges of 
this parameter. Thus we expect chaotic behaviour as well for values of $A$ in 
certain intervals.
For the area parameter with value $A=0.55$, the equilibrium point of the 
system \eqref{redsistp} is stable. If this equilibrium point is unique, then 
the 
system must have a hidden strange attractor. Examples of this type of dynamical 
systems are rather limited (cf. \cite{MoJaSpGo2013}). Further analysis in this 
direction as well as for the case $A=0.65$, and the implications of 
the results in this paper to the study of ``cavitation'' in fluids mentioned in 
the introduction, shall be pursued elsewhere.

\noindent\textbf{Acknowledgements:}
This research was sponsored in part by the NSF--PREM Program of the UPRH (Grant No. DMR--1523463).

\end{document}